\documentclass[2p]{article}
\usepackage{amssymb,amsmath,amsthm}
\usepackage{indentfirst}
\usepackage{exscale}
\usepackage{relsize}
\usepackage{xcolor}
\usepackage[numbers,sort&compress]{natbib}

\usepackage{geometry}
\newcommand{\R}{\mathbb{R}}
\theoremstyle{definition}

\allowdisplaybreaks
\theoremstyle{remark}

\numberwithin{equation}{section}
 \UseRawInputEncoding

\begin{document}
 \title{\Large\bf{ Infinitely many  solutions for a class of elliptic boundary value problems with $(p,q)$-Kirchhoff type }
 }
 \date{}
\author { \ Zongxi Li$^{1}$,  Wanting Qi$^{1}$, Xingyong Zhang$^{1,2}$\footnote{Corresponding author, E-mail address: zhangxingyong1@163.com}\\\
      {\footnotesize $^{1}$Faculty of Science, Kunming University of Science and Technology, Kunming, Yunnan, 650500, P.R. China.}\\
      {\footnotesize $^{1,2}$Research Center for Mathematics and Interdisciplinary Sciences, Kunming University of Science and Technology,}\\
 {\footnotesize Kunming, Yunnan, 650500, P.R. China.}\\
 }
 \date{}
 \maketitle
 \begin{center}
 \begin{minipage}{15cm}
 \par
 \small  {\bf Abstract:} In this paper,
 we investigate the existence of infinitely many solutions for
  the following  elliptic boundary value problem with
  $(p,q)$-Kirchhoff type
 \begin{equation*}
\label{a1} \begin{cases}
   -\Big[M_1\left(\int_\Omega|\nabla u_1|^p dx\right)\Big]^{p-1}\Delta_p u_1+\Big[M_3\left(\int_\Omega a_1(x)|u_1|^p dx\right)\Big]^{p-1}a_1(x)|u_1|^{p-2}u_1=G_{u_1}(x,u_1,u_2)\ \ \mbox{in }\Omega,\\
   -\Big[M_2\left(\int_\Omega|\nabla u_2|^q dx\right)\Big]^{q-1}\Delta_q u_2+\Big[M_4\left(\int_\Omega a_2(x)|u_2|^q dx\right)\Big]^{q-1}a_2(x)|u_2|^{q-2}u_2=G_{u_2}(x,u_1,u_2)\ \ \mbox{in }\Omega,\\
   u_1=u_2=0\ \ \quad \quad \quad \quad \quad \quad \quad \  \mbox{ on   }\partial\Omega.
   \end{cases}
 \end{equation*}
 By using a critical point theorem due to Ding in
 [Y. H. Ding, Existence and multiplicity results for homoclinic solutions to a class of Hamiltonian systems.
 Nonlinear Anal, 25(11)(1995)1095-1113], we obtain that system has
 infinitely many solutions under the sub-$(p,q)$ conditions.

 \end{minipage}
 \end{center}

  \vskip2mm
 {\section{Introduction and Main results }}
 \allowdisplaybreaks
  In this paper,
  we investigate the following nonlocal elliptic boundary value problem with
  $(p,q)$-Kirchhoff type
\begin{equation}
\label{a1} \begin{cases}
   -\Big[M_1\left(\int_\Omega|\nabla u_1|^p dx\right)\Big]^{p-1}\Delta_p u_1+\Big[M_3\left(\int_\Omega a_1(x)|u_1|^p dx\right)\Big]^{p-1}a_1(x)|u_1|^{p-2}u_1=G_{u_1}(x,u_1,u_2)\ \ \mbox{in }\Omega,\\
   -\Big[M_2\left(\int_\Omega|\nabla u_2|^q dx\right)\Big]^{q-1}\Delta_q u_2+\Big[M_4\left(\int_\Omega a_2(x)|u_2|^q dx\right)\Big]^{q-1}a_2(x)|u_2|^{q-2}u_2=G_{u_2}(x,u_1,u_2)\ \ \mbox{in }\Omega,\\
   u_1=u_2=0\ \ \quad \quad \quad \quad \quad \quad \quad \  \mbox{ on   }\partial\Omega,
\end{cases}
 \end{equation}
 where $1<p\le \frac{pN}{N-p} ,1<q\le \frac{qN}{N-q}$, $N$ is a integer with $N>\max\{p,q\}$, $\Delta_pu_1=\mbox{div}(|\nabla u_1|^{p-2}\nabla u_1)$,  $\Delta_qu_2=\mbox{div}(|\nabla u_2|^{q-2}\nabla u_2)$, $\Omega$ is a bounded domain of $\mathbb R^N$ with smooth
 boundary $\partial \Omega$, $M_i:\mathbb R^+\to\mathbb R$, $i=1,2,3,4$ are continuous functions and
 $a_i\in C(\bar{\Omega},\mathbb R^+)$, $i=1,2$.

 \vskip2mm
 \par
 In 2011, Cheng, Wu and Liu \cite{BTC} investigated the following nonlocal elliptic
 system of  $(p,q)$-Kirchhoff type with a parameter $\lambda$:
 \begin{equation}
\label{a2} \begin{cases}
   -\Big[M_1\left(\int_\Omega|\nabla u|^p dx\right)\Big]^{p-1}\Delta_p u=\lambda F_{u}(x,u,v)\ \ \mbox{in }\Omega,\\
   -\Big[M_2\left(\int_\Omega|\nabla v|^q dx\right)\Big]^{q-1}\Delta_q v=\lambda  F_{v}(x,u,v)\ \ \mbox{in }\Omega,\\
   u=v=0\ \ \quad \quad \quad \quad \quad \quad \quad \  \mbox{ on   }\partial\Omega,
   \end{cases}
 \end{equation}
 where $\Omega\subset \mathbb R^N(N\ge 1)$ is a bounded smooth domain, $\lambda \in
 (0,+\infty)$, $p>N$, $q>N$, $M_i:\mathbb R^+\to\mathbb R^+$,
 $i=1,2$ are continuous functions with bounded conditions.
 Under some reasonable conditions, by using a critical point theorem due to Bonanno in \cite{Bonanno}, they obtained
 system (\ref{a2}) has at least two weak solutions, and by using an equivalent formulation (\cite{Bonanno2}, Theorem 2.3) of a three critical points theorem due to Ricceri in \cite{Ricceri}, they obtained system (\ref{a2}) has at least three weak solutions.
 Subsequently,
  Chen et al. \cite{ChenGS} and Massar et al. \cite{Massar2015} both investigated the following the following nonlocal elliptic system of  $(p,q)$-Kirchhoff type with two parameters $\lambda$ and $\mu$:
  \begin{equation}
\label{a3} \begin{cases}
   -\Big[M_1\left(\int_\Omega|\nabla u|^p dx\right)\Big]^{p-1}\Delta_p u=\lambda F_{u}(x,u,v)+\mu G_{u}(x,u,v)\ \ \mbox{in }\Omega,\\
   -\Big[M_2\left(\int_\Omega|\nabla v|^q dx\right)\Big]^{q-1}\Delta_q v=\lambda  F_{v}(x,u,v)+\mu  G_{v}(x,u,v)\ \ \mbox{in }\Omega,\\
   u=v=0\ \ \quad \quad \quad \quad \quad \quad \quad \  \mbox{ on   }\partial\Omega.
   \end{cases}
 \end{equation}
 Under different conditions for $F$ and $G$, by using a critical point theorem due to Ricceri in
 \cite{Ricceri2}, they obtained that system (\ref{a3}) has at least three weak solutions and they generalized the corresponding result in \cite{BTC}.
 Then, the results of \cite{ChenGS} were extended slightly to a Dirichlet boundary problem involving the $(p_{1}, ..., p_{n})$-Kirchhoff type systems in \cite{Fang2014}.

 \par
 Moreover, in \cite{Chung}, Chung investigated the following system with a parameter $\lambda$:
  \begin{equation}
\label{a4} \begin{cases}
   -M_1\left(\int_\Omega|\nabla u|^p dx\right)\Delta_p u=\lambda a(x)f(u,v)\ \ \mbox{in }\Omega,\\
   -M_2\left(\int_\Omega|\nabla v|^q dx\right)\Delta_q v=\lambda b(x)g(u,v)\ \ \mbox{in }\Omega,\\
   u=v=0\ \ \quad \quad \quad \quad \quad \quad \quad \  \mbox{ on   }\partial\Omega,
   \end{cases}
 \end{equation}
where $a,b\in C(\bar{\Omega})$.
 By using sub and supersolutions method, under some reasonable conditions for $f$ and $g$, the author obtained that system
 (\ref{a4}) has a positive solutions when $\lambda>\lambda^*$ for some $\lambda^*>0$.
 For problem (\ref{a4}), there are many similar researches and further results, we refer to the literature \cite{Shakeri2021,Rasouli2016} and  references therein.

 \par
 As pointed out in \cite{BTC}, (\ref{a2}) is related to the Kirchhoff
 equation:
 \begin{equation}
 \rho\frac{\partial^2 u}{\partial t^2}-\left(\frac{P_0}{h}+\frac{E}{2L}\int_0^L\left|\frac{\partial u}{\partial
 x}\right|^2dx\right)\frac{\partial^2u}{\partial x^2}=0,
 \end{equation}
 where $\rho$, $P_0$, $h$, $E$ and $L$ are parameters with special
 meanings. Kirchhoff problems have been studied extensively
 and lots of interesting results have been obtained by using variational
 method. We refer to \cite{Alves}, \cite{Bahrouni}, \cite{Bensedik}, \cite{ChenC}, \cite{Correa}, \cite{MaTF} and \cite{SunJT}.

 \par
 In this paper, motivated by \cite{BTC}, \cite{ChenGS} and \cite{Massar2015}, we
 investigate the existence of infinitely many solutions for system
 (\ref{a1}). Via a critical point theorem due to Ding in \cite{Ding1995}, we establish two
 results under sub-$(p,q)$ conditions.  To be
 precise, we obtain the following results:
 \vskip1mm
 \noindent
 {\bf Theorem 1.1.}\ \   {\it Assume that the following
 conditions hold:\\
 ($\mathcal{M}$)\ \ there exist positive constants $\theta\in(0,\min\{p,q\})$, $m_*$ and $m^*$  such that
 \begin{eqnarray*}
         \frac{m_*}{m^*}
  \ge  \max\left\{\left(\frac{\theta\max\limits_{x\in\bar{\Omega}}a_1(x)}{p\min\limits_{x\in{\bar{\Omega}}}a_1(x)}\right)^{\frac{1}{p-1}},
       \left(\frac{\theta\max\limits_{x\in\bar{\Omega}}a_2(x)}{q\min\limits_{x\in{\bar{\Omega}}}a_2(x)}\right)^{\frac{1}{q-1}}\right\}
 \end{eqnarray*}
 and
  \begin{eqnarray*}
 　m_*\le M_i(s)\le m^*, \forall s\ge 0,i=1,2,3,4;
 \end{eqnarray*}
 ($\mathcal{A}$)\ \ $a_i\in C(\bar{\Omega},\mathbb R^+)$ and $\inf\limits_{\Omega}a_i(x)> 0$, $i=1,2$;\\
 (G0) \ \  $G(x,0,0)\equiv 0$ and  $G(x,z_1,z_2)=G(x,-z_1,-z_2)$;\\
 (G1)\ \ $G(x,z_1,z_2)\in C^1(\bar{\Omega}\times\mathbb R\times\mathbb R,\mathbb R)$ and there exist constants $c_1, c_2>0$, $s_1\in \left[p,\frac{(p-1)N}{N-p}\right)$ and $s_2\in\left[q,\frac{(q-1)N}{N-q}\right)$
 such that
 \begin{eqnarray*}
  & &  |G_{z_1}(x,z_1,z_2)|\le c_1(1+|z_1|^{s_1}), \quad \mbox{for all } x\in \Omega\times\mathbb R\times\mathbb R,\\
  & &  |G_{z_2}(x,z_1,z_2)|\le c_2(1+|z_2|^{s_2}), \quad \mbox{for all } x\in \Omega\times\mathbb R\times\mathbb R;
 \end{eqnarray*}
 (G2)\ \ $$
           \lim_{|z_1|+|z_2|\to\infty}\frac{G(x,z_1,z_2)}{|z_1|^p+|z_2|^q}< \min\left\{\frac{m_*^{p-1}\inf\limits_{x\in \Omega}a_1(x)}{p},\frac{m_*^{q-1}\inf\limits_{x\in \Omega}a_2(x)}{q}\right\}\  \   \mbox{uniformly for  } x\in
           \Omega;
           $$
 (G3)\ \ there exist $\gamma_1\in [1,p)$, $\gamma_2\in [1,q)$ and $C^*>0$ such that
 $$
  G(x,z_1,z_2)\ge C^*(|z_1|^{\gamma_1}+|z_2|^{\gamma_2}),\quad \mbox{for all }(x,z_1,z_2)\in \Omega\times  \mathbb R\times \mathbb R;
 $$
 (G4)\ \ there exists a function $h\in L^1(\Omega,\R)$ such that
 $$
   \theta G(x,z_1,z_2)-G_{z_1}(x,z_1,z_2)z_1-G_{z_2}(x,z_1,z_2)z_2\ge h(t) \ \ \ \ \mbox{for a.e. }  x\in
  \Omega,
 $$
 and
 $$\lim_{|z_1|+|z_2|\to\infty}[\theta G(x,z_1,z_2)-G_{z_1}(x,z_1,z_2)z_1-G_{z_2}(x,z_1,z_2)z_2]=+\infty,\ \  \mbox{for  a.e. } x\in \Omega.$$
  Then system (\ref{a1}) has  infinitely many nontrivial solutions.  }

 \vskip1mm
 \noindent
 {\bf Theorem 1.2.}\ \   {\it Assume that ($\mathcal{M}$), ($\mathcal{A}$), (G0), (G1), (G3), (G4)
 and the following condition hold:\\
 (G2)$'$\ \ $$
           \lim_{|z_1|+|z_2|\to\infty}\frac{G(x,z_1,z_2)}{|z_1|^p+|z_2|^q}<\min\left\{\frac{\min\{m_*, m_*^{p-1}\inf\limits_{x\in \Omega}a_1(x)\}}{p\tau_{p,p}},\frac{\min\{m_*, m_*^{q-1}\inf\limits_{x\in \Omega}a_2(x)\}}{q\tau_{q,q}}\right\}
           $$
  uniformly for  $x\in \Omega$, where $\tau_{p,p}$ and  $\tau_{q,q}$ are  the  embedding constants in $W_0^{1,p}(\Omega)\hookrightarrow L^p(\Omega)$ and  $W_0^{1,q}(\Omega)\hookrightarrow L^q(\Omega)$, respectively. Then system (\ref{a1}) has  infinitely many nontrivial solutions.  }

 {\section{Preliminaries}}
  \setcounter{equation}{0}
  \par
 Let $r>1$. On  $W_0^{1,r}(\Omega)$, define the norm
   $$
   \|u\|_{1,r}=\left(\int_\Omega |\nabla u|^r dx+\int_\Omega | u|^r dx\right)^{1/r}.
  $$
  Then $(W_0^{1,r}(\Omega),\|\cdot\|_{1,r})$ is a reflexive and separable  Banach
  space. It is well known that there exist $\{v_n\}_{n\in\mathbb N}\subset W_0^{1,r}(\Omega) $ such
  that $\overline{\mbox{span}}\{v_n;n\in\mathbb
  N\}=W_0^{1,r}(\Omega)$. Let $X_{r,j}=\mathbb R v_j$. Then $W_0^{1,r}(\Omega)=\oplus_{j\ge 1}X_{r,j}$.
  Define
  \begin{eqnarray}
   E_{r,k}^{(1)}=\oplus_{j=1}^k X_{r,j},\quad  E_{r,k}^{(2)}=\overline{\oplus_{j\ge k}
   X_{r,j}}.
  \end{eqnarray}
  Then $W_0^{1,r}(\Omega)=E_{r,k}^{(1)}\oplus E_{r,k}^{(2)}$ (see \cite{sbl}).
  \par
   Let $\mathcal{W}=W_0^{1,p}(\Omega)\times W_0^{1,q}(\Omega)$. Then
   $$
   \mathcal{W}= \left(E_{p,1}^{(1)}\oplus E_{p,1}^{(2)}\right)\times \left(E_{q,1}^{(1)}\oplus E_{q,1}^{(2)}\right)
   =\left(E_{p,1}^{(1)}\times E_{q,1}^{(1)}\right)\oplus \left( E_{p,1}^{(2)}\times   E_{q,1}^{(2)}\right).
   $$
   On $\mathcal{W}$,  define
  the norm
  $$
   \|(u_1,u_2)\|=\|u_1\|_{1,p}+\|u_2\|_{1,q},
  $$
  where $u_1\in W_0^{1,p}(\Omega)$ and  $u_2\in  W_0^{1,q}(\Omega)$.
  Then $(\mathcal{W}, \|\cdot\|)$ is also a reflexive Banach  space.

   \vskip2mm
 \noindent
 {\bf Lemma 2.1.} (see  \cite{Adams})\ \  {\it $W_0^{1,p}$ ($W_0^{1,q}$) is compactly
 embedded in $L^r\equiv L^r(\Omega)$ for $p\le r<\frac{pN}{N-p}$ ($q\le r<\frac{qN}{N-q}$) and
 continuously embedded in $L^{\frac{pN}{N-p}}\equiv
 L^{\frac{pN}{N-p}}(\Omega)$ $\left(L^{\frac{qN}{N-q}}\equiv
 L^{\frac{qN}{N-q}}(\Omega)\right)$, and hence for every $p\le r\le
 \frac{pN}{N-p}$ $\left(q\le r\le
 \frac{qN}{N-q}\right)$, there exists $\tau_{r,p}>0$ $( \tau_{r,q}>0 )$ such that
 $$
  |u|_{r}\le \tau_{r,p}\|u\|_{1,p}\quad (|u|_r\le \tau_{r,q}\|u\|_{1,q}), \quad \forall u\in W_0^{1,p} (u\in W_0^{1,q}),
 $$
 where $|\cdot|_r$ denotes the usual norm in $L^r$ for all $p\le r\le
 \frac{pN}{N-p}$ $\left(q\le r\le
 \frac{qN}{N-q}\right)$.}

  \par
  Define
  \begin{eqnarray*}
  &     &  \hat{M}_1(t)=\int_0^t[M_1(s)]^{p-1}ds,\quad  \hat{M}_2(t)=\int_0^t[M_2(s)]^{q-1}ds,\\
  &     &  \hat{M}_3(t)=\int_0^t[M_3(s)]^{p-1}ds,\quad  \hat{M}_4(t)=\int_0^t[M_4(s)]^{q-1}ds,\quad \forall t\ge  0.
  \end{eqnarray*}
  On $\mathcal{W}$, define the functional $I:\mathcal{W}\to\mathbb R$ by
  \begin{eqnarray}\label{b1}
          I(u)
 &  =  &  I(u_1,u_2) \nonumber\\
 &  =  &  \varphi(u)+\psi(u) \nonumber\\
 &  =  &  \varphi(u_1,u_2)+\psi(u_1,u_2)\nonumber\\
 &  =  &  -\frac{1}{p}\hat{M}_1\left(\int_\Omega |\nabla u_1|^pdx\right)-\frac{1}{q}\hat{M}_2\left(\int_\Omega |\nabla u_2|^qdx\right) \nonumber\\
 &     &　　- \frac{1}{p}\hat{M}_3\left(\int_\Omega a_1(x)|u_1|^p dx\right)-\frac{1}{q}\hat{M}_4\left( \int_\Omega a_2(x)|u_2|^q dx\right)\nonumber\\
 &     &  +\int_{\Omega}G(x,u_1,u_2)dx,\quad \forall u=(u_1,u_2)\in \mathcal{W},
  \end{eqnarray}
 where
 \begin{eqnarray*}
           \varphi(u)=  \varphi(u_1,u_2)
   &  =  & -\frac{1}{p}\hat{M}_1\left(\int_\Omega |\nabla u_1|^pdx\right)-\frac{1}{q}\hat{M}_2\left(\int_\Omega |\nabla u_2|^qdx\right) \nonumber\\
   &     &　　- \frac{1}{p}\hat{M}_3\left(\int_\Omega a_1(x)|u_1|^p dx\right)-\frac{1}{q}\hat{M}_4\left( \int_\Omega a_2(x)|u_2|^q dx\right),\nonumber\\
  \psi(u)=\psi(u_1,u_2)&  =  & \int_\Omega G(x,u_1,u_2)dx,\quad \forall u=(u_1,u_2)\in \mathcal{W}.
 \end{eqnarray*}

 \vskip2mm
 \noindent
 {\bf Lemma 2.2.} {\it Suppose that the following condition holds:\\
 (G1)$'$\ \ $G(x,z_1,z_2)\in C^1(\bar{\Omega}\times\mathbb R\times\mathbb R,\mathbb R)$ and there exist constants $c_1, c_2>0$, $s_1\in \left[p,\frac{(p-1)N+p}{N-p}\right)$ and $s_2\in\left[q,\frac{(q-1)N+q}{N-q}\right)$
 such that
 \begin{eqnarray*}
  & &  |G_{z_1}(x,z_1,z_2)|\le c_1(1+|z_1|^{s_1}), \quad \mbox{for all } x\in \Omega\times\mathbb R\times\mathbb R,\\
  & &  |G_{z_2}(x,z_1,z_2)|\le c_2(1+|z_2|^{s_2}), \quad \mbox{for all } x\in \Omega\times\mathbb R\times\mathbb  R.
 \end{eqnarray*}
  Then $\psi\in C^1(\mathcal{W},\R)$ and $\psi':\mathcal{W}\to \mathcal{W}^*$ is compact, and hence $\varphi\in C^1(\mathcal{W},\R)$. Moreover,
   \begin{eqnarray*}
    \langle \varphi'(u_1,u_2), (v_1,v_2)\rangle
    &  = & -\Big[M_1\left(\int_\Omega|\nabla u_1|^p dx\right)\Big]^{p-1}\int_\Omega|\nabla u_1|^{p-2}(\nabla    u_1,\nabla    v_1)dx\\
    &    &    -\Big[M_2\left(\int_\Omega|\nabla u_2|^q dx\right)\Big]^{q-1}\int_\Omega|\nabla u_2|^{q-2}(\nabla    u_2,\nabla v_2)dx  \\
    &    &    -\Big[M_3\left(\int_\Omega a_1(x)|u_1|^p dx\right)\Big]^{p-1}\int_\Omega a_1(x)|u_1|^{p-2}(u_1,v_1)dx \\
    &    &    -\Big[M_4\left(\int_\Omega a_2(x)|u_2|^q dx\right)\Big]^{q-1}\int_\Omega a_2(x)|u_2|^{q-2}(u_2,v_2)dx
  \end{eqnarray*}
  and
  \begin{eqnarray*}
    \langle\psi'(u_1,u_2), (v_1,v_2)\rangle
    &  = & \int_\Omega G_{z_1}(x,u_1,u_2)v_1dx+\int_\Omega G_{z_2}(x,u_1,u_2)v_2dx
   \end{eqnarray*}
  and
   \begin{eqnarray} \label{dd5}
    \langle I'(u_1,u_2), (v_1,v_2)\rangle      =  \langle\varphi'(u_1,u_2), (v_1,v_2)\rangle +\langle\psi'(u_1,u_2), (v_1,v_2)\rangle
  \end{eqnarray}
  for all $(u_1,u_2),(v_1,v_2)\in \mathcal{W}$ and critical points of $I$ are solutions
  of system (\ref{a1}).}
 \vskip2mm
 \noindent
 {\bf Proof}. The proof is essentially identical to that of Proposition B.10 in \cite{Ra}. Hence we omit the details.

 \vskip2mm
 \noindent
 {\bf Remark 2.2.} Obviously, (G1) implies that (G1)$'$. Hence,
 Lemma 2.2 also holds under (G1).

 \vskip2mm
 \noindent
 {\bf Lemma 2.3.} (see  \cite{Ding1995}, Lemma 2.4)\ \ {\it Let $E$ be an infinite
 dimensional Banach space and let $f\in C^1(E,\R)$ be even, satisfy
 (PS), and $f(0)=0$. If $E=E_1\oplus E_2$, where $E_1$ is finite
 dimensional, and $f$ satisfies
 \par
 ($f_1$)\ \ $f$ is bounded from above on $E_2$,
 \par
 ($f_2$)\ \ for each finite dimensional subspace $\tilde{E}\subset
 E$, there are positive constants $\rho=\rho(\tilde{E})$ and
 $\sigma=\sigma(\tilde{E})$ such that $f\ge 0$ on $B_\rho\cap
 \tilde{E}$ and $f|_{\partial B_\rho\cap \tilde{E}}\ge \sigma$ where
 $B_\rho=\{x\in E;\|x\|\le \rho\}$, then $f$ possesses infinitely
 many nontrivial critical points.}

  \vskip2mm
 \noindent
 {\bf Remark 2.1.} As shown in \cite{BBF}, a deformation lemma can be proved with
 replacing the usual (PS)-condition with the (C)-condition introduced by Cerami in \cite{CE}, and it turns out that Lemma 2.2 are true under the (C)-condition. We say that $\varphi$ satisfies the (C)-condition, i.e. for every sequence $\{u_n\}\subset E$, $\{u_n\}$ has
 a convergent subsequence if $\varphi(u_n)$ is bounded and $(1 + \|u_n\|)\|\varphi'(u_n)\|\to 0$ as $n \to \infty.$

 \vskip4mm
 \noindent
  {\section{Proof}}
  \setcounter{equation}{0}
  \vskip2mm
 \noindent
 {\bf Lemma 3.1.} {\it Assume that (G2) holds. Then $I$ is bounded from above on $\mathcal{W}$.}
 \vskip1mm
 \noindent
 {\bf Proof.} It follows from (G2) that there exist positive
 constants
 $$
 \varepsilon< \min\left\{\frac{m_*^{p-1}\inf\limits_{x\in \Omega}a_1(x)}{p},\frac{m_*^{q-1}\inf\limits_{x\in \Omega}a_2(x)}{q}\right\}
 $$
 and $r_0>0$ such that
  \begin{eqnarray}\label{b3}
   G(x,z_1,z_2)\le \varepsilon (|z_1|^p+|z_2|^q), \quad \mbox{for all }|z_1|+|z_2|>r_0 \mbox{ and all }x\in \Omega.
   \end{eqnarray}
  Then (\ref{b3}) and the continuity of $G$ imply that there exists
  a positve constant $C_0$ such that
 \begin{eqnarray}\label{b4}
   G(x,z_1,z_2)\le \varepsilon (|z_1|^p+|z_2|^q)+C_0, \quad \mbox{for all } (x,z_1,z_2)\in \Omega\times \mathbb R\times \mathbb R.
   \end{eqnarray}
 Hence, by ($\mathcal{M}$), ($\mathcal{A}$), (\ref{b4}) and Lemma 2.1, for all $u\in \mathcal{W}$, we have
\begin{eqnarray}\label{b5}
          I(u)
 &  =  &  I(u_1,u_2) \nonumber\\
 &  =  &  \varphi(u_1,u_2)+\psi(u_1,u_2)\nonumber\\
 &  =  &  -\frac{1}{p}\hat{M}_1\left(\int_\Omega |\nabla u_1|^pdx\right)-\frac{1}{q}\hat{M}_2\left(\int_\Omega |\nabla u_2|^qdx\right) \nonumber\\
 &     &　　- \frac{1}{p}\hat{M}_3\left(\int_\Omega a_1(x)|u_1|^p dx\right)-\frac{1}{q}\hat{M}_4\left( \int_\Omega a_2(x)|u_2|^q dx\right)\nonumber\\
 &     &  +\int_{\Omega}G(x,u_1,u_2)dx\nonumber\\
 & \le &  -\frac{m_*^{p-1}}{p}\int_\Omega |\nabla u_1|^pdx-\frac{m_*^{q-1}}{q}\int_\Omega |\nabla u_2|^qdx \nonumber\\
 &     &　　- \frac{m_*^{p-1}}{p}\int_\Omega a_1(x)|u_1|^p dx-\frac{m_*^{q-1}}{q} \int_\Omega a_2(x)|u_2|^q dx\nonumber\\
 &     &  +\varepsilon\int_{\Omega}|u_1|^pdx+\varepsilon\int_{\Omega}|u_2|^qdx+C_0\Omega\nonumber\\
 & \le &  - \frac{m_*^{p-1}\inf\limits_{x\in \Omega}a_1(x)}{p}\int_\Omega |u_1|^p dx-\frac{m_*^{q-1}\inf\limits_{x\in \Omega}a_2(x)}{q} \int\limits_\Omega |u_2|^q dx\nonumber\\
 &     & +\varepsilon\int_{\Omega}|u_1|^pdx+\varepsilon\int_{\Omega}|u_2|^qdx+C_0\Omega.
  \end{eqnarray}
 Note that  $\varepsilon<\min\left\{\frac{m_*^{p-1}\inf\limits_{x\in \Omega}a_1(x)}{p},\frac{m_*^{q-1}\inf\limits_{x\in \Omega}a_2(x)}{q}\right\}$. Hence, (\ref{b5})
 implies the conclusion holds.

  \vskip2mm
 \noindent
 {\bf Lemma 3.2.} {\it Assume that (G2)$'$ holds. Then $I$ is bounded from above on $\mathcal{W}$.}
 \vskip1mm
 \noindent
 {\bf Proof.} It follows from (G2)$'$ that there exist positive
 constants
 $$
 \varepsilon'<\min\left\{\frac{\min\{m_*, m_*^{p-1}\inf\limits_{x\in \Omega}a_1(x)\}}{p\tau_{p,p}},\frac{\min\{m_*, m_*^{q-1}\inf\limits_{x\in \Omega}a_2(x)\}}{q\tau_{q,q}}\right\}
 $$
 and $r_0'>0$ such that
  \begin{eqnarray}\label{b31}
   G(x,z_1,z_2)\le \varepsilon' (|z_1|^p+|z_2|^q), \quad \mbox{for all }|z_1|+|z_2|>r_0' \mbox{ and all }x\in \Omega.
   \end{eqnarray}
  Then (\ref{b31}) and the continuity of $G$ imply that there exists
  a positve constant $C_0'$ such that
 \begin{eqnarray}\label{b41}
   G(x,z_1,z_2)\le \varepsilon' (|z_1|^p+|z_2|^q)+C_0', \quad \mbox{for all } (x,z_1,z_2)\in \Omega\times \mathbb R\times \mathbb R.
   \end{eqnarray}
 Hence, by ($\mathcal{M}$), ($\mathcal{A}$), (\ref{b41}) and Lemma 2.1, for all $u\in \mathcal{W}$, we have
\begin{eqnarray}\label{b51}
          I(u)
 &  =  &  -\frac{1}{p}\hat{M}_1\left(\int_\Omega |\nabla u_1|^pdx\right)-\frac{1}{q}\hat{M}_2\left(\int_\Omega |\nabla u_2|^qdx\right) \nonumber\\
 &     &　　- \frac{1}{p}\hat{M}_3\left(\int_\Omega a_1(x)|u_1|^p dx\right)-\frac{1}{q}\hat{M}_4\left( \int_\Omega a_2(x)|u_2|^q dx\right)\nonumber\\
 &     &  +\int_{\Omega}G(x,u_1,u_2)dx\nonumber\\
 & \le &  -\frac{m_*^{p-1}}{p}\int_\Omega |\nabla u_1|^pdx-\frac{m_*^{q-1}}{q}\int_\Omega |\nabla u_2|^qdx \nonumber\\
 &     &　　- \frac{m_*^{p-1}}{p}\int_\Omega a_1(x)|u_1|^p dx-\frac{m_*^{q-1}}{q} \int_\Omega a_2(x)|u_2|^q dx\nonumber\\
 &     &  +\varepsilon'\int_{\Omega}|u_1|^pdx+\varepsilon'\int_{\Omega}|u_2|^qdx+C_0'\Omega\nonumber\\
 & \le &  - \min\left\{\frac{m_*^{p-1}}{p},\frac{m_*^{p-1}\inf\limits_{x\in \Omega}a_1(x)}{p}\right\}\|u_1\|_{1,p}^p-\min\left\{\frac{m_*^{q-1}}{q},\frac{m_*^{q-1}\inf\limits_{x\in \Omega}a_2(x)}{q} \right\}\|u_2\|_{1,q}^q\nonumber\\
 &     & +\varepsilon'\tau_{p,p}^p\|u_1\|_{1,p}^p+\varepsilon'\tau_{q,q}^q\|u_2\|_{1,q}^q+C_0'\Omega.
  \end{eqnarray}
 Note that  $\varepsilon'<\min\left\{\frac{\min\{m_*, m_*^{p-1}\inf\limits_{x\in \Omega}a_1(x)\}}{p\tau_{p,p}},\frac{\min\{m_*, m_*^{q-1}\inf\limits_{x\in \Omega}a_2(x)\}}{q\tau_{q,q}}\right\}$. Hence, (\ref{b51})
 implies the conclusion holds.

  \vskip2mm
 \noindent
 {\bf Lemma 3.3.} {\it  Assume that (G3) holds. Then for each finite dimensional subspace $\tilde{\mathcal{W}}\subset
 \mathcal{W}$, there are positive constants $\rho=\rho(\tilde{\mathcal{W}})$ and
 $\sigma=\sigma(\tilde{\mathcal{W}})$ such that $I\ge 0$ on $B_\rho\cap
 \tilde{\mathcal{W}}$ and $I|_{\partial B_\rho\cap \tilde{\mathcal{W}}}\ge \sigma$.}

 \vskip1mm
 \noindent
 {\bf Proof.}
 Since $\tilde{\mathcal{W}}$ is finite dimensional, all norms on $\tilde{\mathcal{W}}$
 are equivalent. Hence there exist $d_i=d_i(\tilde{\mathcal{W}})>0,i=1,2,3,4$ such
 that for all $u=(u_1,u_2)\in \tilde{\mathcal{W}}$,
  \begin{eqnarray}\label{b7}
 &  &　d_1\|u_1\|_{1,p}\le\left(\int_\Omega |u_1|^{\gamma_1} dx\right)^{1/\gamma_1}\le d_2\|u_1\|_{1,p},\nonumber\\
 &  &　d_3\|u_2\|_{1,q}\le\left(\int_\Omega |u_2|^{\gamma_2} dx\right)^{1/\gamma_2}\le d_4\|u_2\|_{1,q}.
  \end{eqnarray}
 It follows from (G3) that
 \begin{eqnarray}\label{b8}
           I(u)
  &  =  &  -\frac{1}{p}\hat{M}_1\left(\int_\Omega |\nabla u_1|^pdx\right)-\frac{1}{q}\hat{M}_2\left(\int_\Omega |\nabla u_2|^qdx\right) \nonumber\\
 &     &　　- \frac{1}{p}\hat{M}_3\left(\int_\Omega a_1(x)|u_1|^p dx\right)-\frac{1}{q}\hat{M}_4\left( \int_\Omega a_2(x)|u_2|^q dx\right)\nonumber\\
 &     &  +\int_{\Omega}G(x,u_1,u_2)dx\nonumber\\
 & \ge & -\frac{(m^*)^{p-1}}{p}\int_\Omega |\nabla u_1|^pdx-\frac{(m^*)^{q-1}}{q}\int_\Omega |\nabla u_2|^qdx \nonumber\\
 &     &　　- \frac{(m^*)^{p-1}}{p}\int_\Omega a_1(x)|u_1|^p dx-\frac{(m^*)^{q-1}}{q} \int_\Omega a_2(x)|u_2|^q dx\nonumber\\
 &     &  +C^{*}\int_{\Omega}|u_1|^{\gamma_1}dx+C^{*}\int_{\Omega}|u_2|^{\gamma_2}dx\nonumber\\
 & \ge & -\frac{(m^*)^{p-1}}{p}\left(1+\max_{x\in \bar{\Omega}}a_1(x)\right)\left(\int_\Omega |\nabla u_1|^pdx+\int_\Omega |u_1|^pdx\right)\nonumber\\
 &     & -\frac{(m^*)^{q-1}}{q}\left(1+\max_{x\in \bar{\Omega}}a_2(x)\right)\left(\int_\Omega |\nabla u_2|^qdx+\int_\Omega |u_2|^qdx\right) \nonumber\\
 &     &  +C^{*}d_1^{\gamma_1}\|u_1\|_{1,p}^{\gamma_1}+C^{*}d_3^{\gamma_2}\|u_2\|_{1,q}^{\gamma_2}\nonumber\\
 & =   & -\frac{(m^*)^{p-1}}{p}\left(1+\max_{x\in \bar{\Omega}}a_1(x)\right)\|u_1\|_{1,p}^p -\frac{(m^*)^{q-1}}{q}\left(1+\max_{x\in \bar{\Omega}}a_2(x)\right)\|u_2\|_{1,q}^q \nonumber\\
 &     &  +C^{*}d_1^{\gamma_1}\|u_1\|_{1,p}^{\gamma_1}+C^{*}d_3^{\gamma_1}\|u_2\|_{1,q}^{\gamma_2}\nonumber\\
 &  =  & \left[C^*d_1^{\gamma_1}-\frac{(m^*)^{p-1}}{p}\left(1+\max_{x\in \bar{\Omega}}a_1(x)\right)\|u_1\|_{1,p}^{p-\gamma_1}\right]\|u_1\|_{1,p}^{\gamma_1}\nonumber\\
 &     & +\left[C^*d_3^{\gamma_2}-\frac{(m^*)^{q-1}}{q}\left(1+\max_{x\in \bar{\Omega}}a_2(x)\right)\|u_2\|_{1,q}^{q-\gamma_2}\right]\|u_2\|_{1,q}^{\gamma_2}.
 \end{eqnarray}
 Let $\rho=\max\left\{\left(\frac{C^*d_1^{\gamma_1}}{(m^*)^{p-1}\left(1+\max\limits_{x\in \bar{\Omega}}a_1(x)\right)}\right)^{1/(p-\gamma_1)},\left(\frac{C^*d_3^{\gamma_2}}{(m^*)^{q-1}\left(1+\max\limits_{x\in \bar{\Omega}}a_2(x)\right)}\right)^{1/(q-\gamma_2)}\right\}$.
 Then $\rho>0$ and there exists a constant $C_1>0$ such that
 $$
  I(u)\ge C_1(\|u_1\|_{1,p}^{\gamma_1}+\|u_2\|_{1,q}^{\gamma_2})\ge 0, \quad \forall \ u\in B_\rho\cap
 \tilde{\mathcal{W}}
 $$
 and it is easy to see that
 $$
  I|_{\partial B_\rho\cap \tilde{{\mathcal{W}}}}>0.
 $$
 Then we complete the proof.

 \vskip2mm
 \noindent
 {\bf Lemma 3.4.} {\it Assume that $(G2)$ and (G4) hold. Then $I$ satisfies the  (C)-condition. }
 \vskip1mm
 \noindent
 {\bf Proof.}    For every sequence $\{u^{[n]}\}_{n\in\mathbb N}=\{(u_1^{[n]},u_2^{[n]})\}_{n\in\mathbb N}\subset \mathcal{W}$, assume that
  $\{I(u^{[n]})\}$ is bounded and  $ \left(1+\|u^{[n]}\|\right)\|I'(u^{[n]})\|\to 0$, as $n\to\infty$.  Then there exists a constant $C_{2}>0$ such that
 \begin{eqnarray}\label{d1}
 |I(u^{[n]})|=|I(u_1^{[n]},u_2^{[n]})|\le C_{2},\ \ \
 \left(1+\|u^{[n]}\|\right)\|I'(u^{[n]})\|\le C_{2},  \
 \ \mbox{for all }n\in \mathbb N.
 \end{eqnarray}
 Then by ($\mathcal{M}$), we have
  \begin{eqnarray}\label{d3}
 &     & (\theta+1)C_2\nonumber\\
 & \ge &  \theta I(u^{[n]})-\langle I'(u^{[n]}),u^{[n]}\rangle\nonumber\\
 &  =  &  -\frac{\theta}{p}\hat{M}_1\left(\int_\Omega |\nabla u_1^{[n]}|^pdx\right)-\frac{\theta}{q}\hat{M}_2\left(\int_\Omega |\nabla u_2^{[n]}|^qdx\right) \nonumber\\
 &     &　　- \frac{\theta}{p}\hat{M}_3\left(\int_\Omega a_1(x)|u_1^{[n]}|^p dx\right)-\frac{\theta}{q}\hat{M}_4\left( \int_\Omega a_2(x)|u_2^{[n]}|^q dx\right)\nonumber\\
 &     &  +\theta\int_{\Omega}G(x,u_1^{[n]},u_2^{[n]})dx\nonumber\\
 &     &  +\Big[M_1\left(\int_\Omega|\nabla u_1^{[n]}|^p dx\right)\Big]^{p-1}\int_\Omega|\nabla u_1^{[n]}|^{p}dx\nonumber\\
 &    &    +\Big[M_2\left(\int_\Omega|\nabla u_2^{[n]}|^q dx\right)\Big]^{q-1}\int_\Omega|\nabla u_2^{[n]}|^{q}dx  \nonumber\\
 &    &    +\Big[M_3\left(\int_\Omega a_1(x)|u_1^{[n]}|^p dx\right)\Big]^{p-1}\int_\Omega a_1(x)|u_1^{[n]}|^{p}dx \nonumber\\
 &    &    +\Big[M_4\left(\int_\Omega a_2(x)|u_2^{[n]}|^q dx\right)\Big]^{q-1}\int_\Omega a_2(x)|u_2^{[n]}|^{q}dx\nonumber\\
 &    &  -\int_\Omega G_{z_1}(x,u_1^{[n]},u_2^{[n]})u_1^{[n]}dx-\int_\Omega G_{z_2}(x,u_1^{[n]},u_2^{[n]})u_2^{[n]}dx\nonumber\\
 & \ge&  -\frac{\theta(m^*)^{p-1}}{p}\int_\Omega |\nabla u_1^{[n]}|^pdx-\frac{\theta(m^*)^{q-1}}{q}\int_\Omega |\nabla u_2^{[n]}|^qdx \nonumber\\
 &     &　　- \frac{\theta(m^*)^{p-1}\max\limits_{x\in\bar{\Omega}}a_1(x)}{p}\int_\Omega |u_1^{[n]}|^p dx\nonumber\\
 &     &  -\frac{\theta(m^*)^{q-1}\max\limits_{x\in\bar{\Omega}}a_2(x)}{q} \int_\Omega |u_2^{[n]}|^q dx\nonumber\\
 &     &  +\theta\int_{\Omega}G(x,u_1^{[n]},u_2^{[n]})dx  +m_*^{p-1}\int_\Omega|\nabla u_1^{[n]}|^{p}dx\nonumber\\
 &    &    +m_*^{q-1}\int_\Omega|\nabla u_2^{[n]}|^{q}dx   +m_*^{p-1}\min_{x\in\bar{\Omega}} a_1(x)\int_\Omega |u_1^{[n]}|^{p}dx \nonumber\\
 &    &    +m_*^{q-1}\min_{x\in\bar{\Omega}} a_2(x)\int_\Omega |u_2^{[n]}|^{q}dx  -\int_\Omega G_{z_1}(x,u_1^{[n]},u_2^{[n]})u_1^{[n]}dx-\int_\Omega G_{z_2}(x,u_1^{[n]},u_2^{[n]})u_2^{[n]}dx\nonumber\\
 & \ge& \theta\int_{\Omega}G(x,u_1^{[n]},u_2^{[n]})dx    -\int_\Omega G_{z_1}(x,u_1^{[n]},u_2^{[n]})u_1^{[n]}dx-\int_\Omega G_{z_2}(x,u_1^{[n]},u_2^{[n]})u_2^{[n]}dx.\nonumber\\
  \end{eqnarray}
 Next we prove that $\{u^{[n]}\}=\{(u_1^{[n]},u_2^{[n]})\}$  is bounded.  Assume that
  \begin{eqnarray}\label{dd1}
  \|u^{[n]}\|=\|u_1^{[n]}\|_{1,p}+\|u_2^{[n]}\|_{1,q}\to \infty, \mbox{ as } n\to  \infty.
  \end{eqnarray}
  Let $z^{[n]}=(z_1^{[n]},z_2^{[n]})=\left(\dfrac{u_1^{[n]}}{\|u_1^{[n]}\|_{1,p}},\dfrac{u_2^{[n]}}{\|u_2^{[n]}\|_{1,q}}\right)$.  Then
  $\|z_1^{[n]}\|_{1,p}= 1$,  $\|z_2^{[n]}\|_{1,q}= 1$ and hence $\|z^{[n]}\|=2$, where $z_1^{[n]}=\dfrac{u_1^{[n]}}{\|u_1^{[n]}\|_{1,p}}$ and $z_2^{[n]}=\dfrac{u_2^{[n]}}{\|u_2^{[n]}\|_{1,q}}$.
 So there exist subsequences, still denoted by $\{z_1^{[n]}\}$ and $\{z_2^{[n]}\}$, such
 that $z_1^{[n]}\rightharpoonup z_1$ and $z_2^{[n]}\rightharpoonup z_2$ on $W_0^{1,p}$ and $W_0^{1,q}$, respectively. Then by Lemma 2.1, we have
  \begin{eqnarray}
 \label{dddd1} &    &  z_1^{[n]}\to z_1 \mbox{ in } L^r \mbox{ for } p\le r<\frac{pN}{N-p} \mbox{ and a.e. } x\in \Omega,\\
 \label{dddd2}&    &  z_2^{[n]}\to z_2 \mbox{ in } L^r \mbox{ for } q\le r<\frac{qN}{N-q}  \mbox{ and a.e. } x\in \Omega.
 \end{eqnarray}
 Thus  by (\ref{b4}), we have
  \begin{eqnarray} \label{c1}
           I(u^{[n]})
    &  =  &  -\frac{1}{p}\hat{M}_1\left(\int_\Omega |\nabla u_1^{[n]}|^pdx\right)-\frac{1}{q}\hat{M}_2\left(\int_\Omega |\nabla u_2^{[n]}|^qdx\right) \nonumber\\
 &     &　　- \frac{1}{p}\hat{M}_3\left(\int_\Omega a_1(x)|u_1^{[n]}|^p dx\right)-\frac{1}{q}\hat{M}_4\left( \int_\Omega a_2(x)|u_2^{[n]}|^q dx\right)  +\int_{\Omega}G(x,u_1^{[n]},u_2^{[n]})dx\nonumber\\
  & \le &   - \min\left\{\frac{m_*^{p-1}}{p},\frac{m_*^{p-1}\inf\limits_{x\in \Omega}a_1(x)}{p}\right\}\|u_1^{[n]}\|_{1,p}^p-\min\left\{\frac{m_*^{q-1}}{q},\frac{m_*^{q-1}\inf\limits_{x\in \Omega}a_2(x)}{q} \right\}\|u_2^{[n]}\|_{1,q}^q\nonumber\\
 &     & +\varepsilon\int_{\Omega}|u_1^{[n]}|^pdx+\varepsilon\int_{\Omega}|u_2^{[n]}|^qdx+C_0\mbox{meas}\Omega\nonumber\\
 & \le &   - \min\left\{\frac{m_*^{p-1}}{p},\frac{m_*^{p-1}\inf\limits_{x\in \Omega}a_1(x)}{p},\frac{m_*^{q-1}}{q},\frac{m_*^{q-1}\inf\limits_{x\in \Omega}a_2(x)}{q} \right\}\left(\|u_1^{[n]}\|_{1,p}^p+\|u_2^{[n]}\|_{1,q}^q\right)\nonumber\\
 &     & +\varepsilon\int_{\Omega}|u_1^{[n]}|^pdx+\varepsilon\int_{\Omega}|u_2^{[n]}|^qdx+C_0\mbox{meas}\Omega.
  \end{eqnarray}
  Then (\ref{c1})  imply that
  \begin{eqnarray}
            \frac{I(u^{[n]})}{\|u_1^{[n]}\|_{1,p}^p+\|u_2^{[n]}\|_{1,q}^q}
 & \le & - \min\left\{\frac{m_*^{p-1}}{p},\frac{m_*^{p-1}\inf\limits_{x\in \Omega}a_1(x)}{p},\frac{m_*^{q-1}}{q},\frac{m_*^{q-1}\inf\limits_{x\in \Omega}a_2(x)}{q} \right\}\nonumber\\
 &     & +\frac{\varepsilon\int_{\Omega}|u_1^{[n]}|^pdx}{\|u_1^{[n]}\|_{1,p}^p+\|u_2^{[n]}\|_{1,q}^q}+\frac{\varepsilon\int_{\Omega}|u_2^{[n]}|^qdx}{\|u_1^{[n]}\|_{1,p}^p+\|u_2^{[n]}\|_{1,q}^q}\nonumber\\
 \label{ddd1}&     & +\frac{C_0\mbox{meas}\Omega}{\|u_1^{[n]}\|_{1,p}^p+\|u_2^{[n]}\|_{1,q}^q}\\
 & \le & - \min\left\{\frac{m_*^{p-1}}{p},\frac{m_*^{p-1}\inf\limits_{x\in \Omega}a_1(x)}{p},\frac{m_*^{q-1}}{q},\frac{m_*^{q-1}\inf\limits_{x\in \Omega}a_2(x)}{q} \right\}\nonumber\\
 \label{dd3}&     & +\frac{\varepsilon\int_{\Omega}|u_1^{[n]}|^pdx}{\|u_1^{[n]}\|_{1,p}^p}+\frac{\varepsilon\int_{\Omega}|u_2^{[n]}|^qdx}{\|u_2^{[n]}\|_{1,q}^q}
  +\frac{C_0\mbox{meas}\Omega}{\|u_1^{[n]}\|_{1,p}^p+\|u_2^{[n]}\|_{1,q}^p}\nonumber\\
 &  =  & - \min\left\{\frac{m_*^{p-1}}{p},\frac{m_*^{p-1}\inf\limits_{x\in \Omega}a_1(x)}{p},\frac{m_*^{q-1}}{q},\frac{m_*^{q-1}\inf\limits_{x\in \Omega}a_2(x)}{q} \right\}\nonumber\\
 &     & +\varepsilon\int_{\Omega}|z_1^{[n]}|^pdx+\varepsilon\int_{\Omega}|z_2^{[n]}|^qdx
  +\frac{C_0\mbox{meas}\Omega}{\|u_1^{[n]}\|_{1,p}^p+\|u_2^{[n]}\|_{1,q}^p}\nonumber\\
 &  =  & - \min\left\{\frac{m_*^{p-1}}{p},\frac{m_*^{p-1}\inf\limits_{x\in \Omega}a_1(x)}{p},\frac{m_*^{q-1}}{q},\frac{m_*^{q-1}\inf\limits_{x\in \Omega}a_2(x)}{q} \right\}\nonumber\\
 &     & +\varepsilon\int_{\Omega}|z_1^{[n]}-z_1+z_1|^pdx+\varepsilon\int_{\Omega}|z_2^{[n]}-z_2+z_2|^qdx
  +\frac{C_0\mbox{meas}\Omega}{\|u_1^{[n]}\|_{1,p}^p+\|u_2^{[n]}\|_{1,q}^p}\nonumber\\
 & \le & - \min\left\{\frac{m_*^{p-1}}{p},\frac{m_*^{p-1}\inf\limits_{x\in \Omega}a_1(x)}{p},\frac{m_*^{q-1}}{q},\frac{m_*^{q-1}\inf\limits_{x\in \Omega}a_2(x)}{q} \right\}\nonumber\\
 &     & +2^{p-1}\varepsilon\int_{\Omega}|z_1^{[n]}-z_1|^pdx+2^{p-1}\varepsilon\int_{\Omega}|z_1|^pdx+2^{q-1}\varepsilon\int_{\Omega}|z_2^{[n]}-z_2|^qdx\nonumber\\
 &     & +2^{q-1}\varepsilon\int_{\Omega}|z_2|^qdx
  +\frac{C_0\mbox{meas}\Omega}{\|u_1^{[n]}\|_{1,p}^p+\|u_2^{[n]}\|_{1,q}^p}.
  \end{eqnarray}
 Let $n\to \infty$. Then by (\ref{d1}), (\ref{dddd1}), (\ref{dddd2}) and  (\ref{dd3}), we get
 \begin{eqnarray}\label{d2}
         \min\left\{\frac{m_*^{p-1}}{p},\frac{m_*^{p-1}\inf\limits_{x\in \Omega}a_1(x)}{p},\frac{m_*^{q-1}}{q},\frac{m_*^{q-1}\inf\limits_{x\in \Omega}a_2(x)}{q} \right\}
    & \le &
          (2^{p-1}+2^{q-1})\varepsilon\left(\int_\Omega|z_1|^pdx+\int_\Omega|z_2|^pdx\right).\nonumber\\
  \end{eqnarray}
 Then it follows from $(\ref{d2})$ that $\int_\Omega|z_1|^pdx+\int_\Omega|z_2|^qdx>0$ and so
 $z=(z_1,z_2)\not=0$.  Next, we claim that the following three conclusions hold:
 \par
 (i)  $\|u_1^{[n]}\|_{1,p}\to\infty$ if $z_1\not= 0$ and $z_2\equiv 0$;
 \par
 (ii) $\|u_2^{[n]}\|_{1,q}\to\infty$ if $z_2\not= 0$ and $z_1\equiv 0$;
 \par
 (iii) $\|u_1^{[n]}\|_{1,p}\to\infty$ or
 $\|u_2^{[n]}\|_{1,q}\to\infty$ if $z_1\not= 0$ and $z_2\not =0$.\\
 Indeed, for conclusion (i), if $\|u_1^{[n]}\|_{1,p}$  is bounded,
 then $\int_\Omega |u_1^{[n]}|^pdx$ is also bounded. Hence, by (\ref{ddd1}) and
 (\ref{dd1}), we have
  \begin{eqnarray*}
            0
 & \le & - \min\left\{\frac{m_*^{p-1}}{p},\frac{m_*^{p-1}\inf\limits_{x\in \Omega}a_1(x)}{p},\frac{m_*^{q-1}}{q},\frac{m_*^{q-1}\inf\limits_{x\in \Omega}a_2(x)}{q} \right\} +\frac{\varepsilon\int_{\Omega}|u_2^{[n]}|^qdx}{\|u_1^{[n]}\|_{1,p}^p+\|u_2^{[n]}\|_{1,q}^q}\nonumber\\
 & \le & - \min\left\{\frac{m_*^{p-1}}{p},\frac{m_*^{p-1}\inf\limits_{x\in \Omega}a_1(x)}{p},\frac{m_*^{q-1}}{q},\frac{m_*^{q-1}\inf\limits_{x\in \Omega}a_2(x)}{q} \right\} +\frac{\varepsilon\int_{\Omega}|u_2^{[n]}|^qdx}{\|u_2^{[n]}\|_{1,q}^q},
  \end{eqnarray*}
  which implies that $\int_{\Omega} |z_2|^qdx>0$ as $n\to\infty$
  that contradicts $z_2\equiv 0$. Similar to the argument of conclusion (i), it is
  easy to obtain conclusion (ii). For conclusion (iii), if $\|u_1^{[n]}\|_{1,p}$ and  $\|u_2^{[n]}\|_{1,p}$  are
  bounded, then   $\int_\Omega |u_1^{[n]}|^pdx$ and $\int_\Omega |u_2^{[n]}|^pdx$ are also
  bounded. Hence, by (\ref{ddd1}) and
 (\ref{dd1}), we have
  \begin{eqnarray*}
            0
 & \le & - \min\left\{\frac{m_*^{p-1}}{p},\frac{m_*^{p-1}\inf\limits_{x\in \Omega}a_1(x)}{p},\frac{m_*^{q-1}}{q},\frac{m_*^{q-1}\inf\limits_{x\in \Omega}a_2(x)}{q} \right\},
  \end{eqnarray*}
 which is a contradiction.
 \par
Let
$${\it{S}}=\{x\in \Omega:\lim\limits_{|z_1|+|z_2|\to\infty}[\theta G(x,z_1,z_2)-G_{z_1}(x,z_1,z_2)z_1-G_{z_2}(x,z_1,z_2)z_2]=+\infty\}$$
and
$${\it{S}_1}={\{x\in \Omega :z(x)=(z_1(x),z_2(x))\not= 0\}}.$$
Then by (G4), we have
 $
  \mbox{meas\;} {\it{S}}>0.
 $
 Then conclusion (i), (ii) and (iii) imply that
  \begin{eqnarray}\label{d6}
  \lim_{n\to\infty}(|u_1^{[n]}(x)|+|u_2^{[n]}(x)|)= \lim_{n\to\infty}(|z_1^{[n]}(x)|\cdot\|u_1^{[n]}\|_{1,p}+|z_2^{[n]}(x)|\cdot\|u_2^{[n]}\|_{1,q})=+\infty\ \ \mbox{for }x\in {\it{S}_1}.
  \end{eqnarray}
  Let
  $J_n(x)=\theta G(x,u_1^{[n]}(x),u_2^{[n]}(x))-G_{z_1}(x,u_1^{[n]}(x),u_2^{[n]}(x))u_1^{[n]}(x)-G_{z_2}(x,u_1^{[n]}(x),u_2^{[n]}(x))u_2^{[n]}(x)$. Then
  (\ref{d6}) and (G4) imply that
  \begin{eqnarray}\label{d7}
   \lim_{n\to \infty} J_n(x)=+\infty \ \ \mbox{for }x\in {\it{S}_1}.
  \end{eqnarray}
 It follows from (\ref{d7}) and Lemma 1 in \cite{TW} that there
 exists a subset ${\it{S}_2}$ of ${\it{S}_1}$ with $\mbox{meas\;} {\it{S}_2}>0$
 such that
 \begin{eqnarray}\label{d8}
   \lim_{n\to \infty} J_n(x)=+\infty \ \ \mbox{uniformly for }x\in {\it{S}_2}.
  \end{eqnarray}
 By (G4), we have
 \begin{eqnarray*}
   &     &  \theta\int_{\Omega}G(x,u_1^{[n]},u_2^{[n]})dx    -\int_\Omega G_{z_1}(x,u_1^{[n]},u_2^{[n]})u_1^{[n]}dx-\int_\Omega G_{z_2}(x,u_1^{[n]},u_2^{[n]})u_2^{[n]}dx\\
   &  =  &   \int_{\it{S}_2} [\theta G(x,u_1^{[n]},u_2^{[n]})-G_{z_1}(x,u_1^{[n]},u_2^{[n]})u_1^{[n]}-G_{z_2}(x,u_1^{[n]},u_2^{[n]})u_2^{[n]}]dx \\
   &     & +\int_{\Omega/\it{S}_2}[\theta G(x,u_1^{[n]},u_2^{[n]})-G_{z_1}(x,u_1^{[n]},u_2^{[n]})u_1^{[n]}-G_{z_2}(x,u_1^{[n]},u_2^{[n]})u_2^{[n]}]dx\\
   & \ge &   \int_{\it{S}_2} [\theta G(x,u_1^{[n]},u_2^{[n]})-G_{z_1}(x,u_1^{[n]},u_2^{[n]})u_1^{[n]}-G_{z_2}(x,u_1^{[n]},u_2^{[n]})u_2^{[n]}]dx
              +\int_{\Omega/\it{S}_2}h(t)dt.
 \end{eqnarray*}
 Let $n\to \infty$. Then by
 Fatou's lemma and (\ref{d8}), we have
 $$
  \int_\Omega[\theta G(x,u_1^{[n]},u_2^{[n]})-G_{z_1}(x,u_1^{[n]},u_2^{[n]})u_1^{[n]}-G_{z_2}(x,u_1^{[n]},u_2^{[n]})u_2^{[n]}]dx\to  +\infty
 $$
 which contradicts (\ref{d3}). Hence $\{u^{[n]}\}=\{(u_1^{[n]},u_2^{[n]})\}$ is bounded and so $\{u_1^{[n]}\}$ and $\{u_2^{[n]}\}$ are also bounded in $W_0^{1,p}$ and  $W_0^{1,q}$, respectively.  Then going
 if necessary to a subsequence, we can assume that
 \begin{eqnarray}\label{dd21}
 u_1^{[n]}\rightharpoonup u_1 \mbox{ in } W_0^{1,p} \mbox{ and } u_2^{[n]}\rightharpoonup u_2 \mbox{ in } W_0^{1,q}.
 \end{eqnarray}
 Then by Lemma 2.1,
 \begin{eqnarray}
 \label{dd4} u_1^{[n]}\to  u_1 \mbox{ in } L^p(\Omega)\quad\mbox{and } u_2^{[n]}\to  u_2 \mbox{ in } L^q(\Omega).
 \end{eqnarray}
 By (\ref{dd5}), we have
 \begin{eqnarray}\label{dd8}
    &    & \langle I'(u_1^{[n]},u_2^{[n]}),    (u_1^{[n]}-u_1,u_2^{[n]}-u_2)\rangle\nonumber\\
    &  = & -\Big[M_1\left(\int_\Omega|\nabla u_1^{[n]}|^p dx\right)\Big]^{p-1}\int_\Omega|\nabla u_1^{[n]}|^{p-2}(\nabla    u_1^{[n]},\nabla  u_1^{[n]}-\nabla u_1)dx\nonumber\\
    &    &    -\Big[M_2\left(\int_\Omega|\nabla u_2^{[n]}|^q dx\right)\Big]^{q-1}\int_\Omega|\nabla u_2^{[n]}|^{q-2}(\nabla    u_2^{[n]},\nabla u_2^{[n]}-\nabla u_2)dx \nonumber\\
    &    &    -\Big[M_3\left(\int_\Omega a_1(x)|u_1^{[n]}|^p dx\right)\Big]^{p-1}\int_\Omega a_1(x)|u_1^{[n]}|^{p-2}(u_1^{[n]},u_1^{[n]}-u_1)dx \nonumber\\
    &    &    -\Big[M_4\left(\int_\Omega a_2(x)|u_2^{[n]}|^q dx\right)\Big]^{q-1}\int_\Omega    a_2(x)|u_2^{[n]}|^{q-2}(u_2^{[n]},u_2^{[n]}-u_2)dx\nonumber\\
    &    & +\int_\Omega G_{z_1}(x,u_1^{[n]},u_2^{[n]})(u_1^{[n]}-u_1)dx+\int_\Omega    G_{z_2}(x,u_1^{[n]},u_2^{[n]})(u_2^{[n]}-u_2)dx.
  \end{eqnarray}
  The boundedness of  $\{u_1^{[n]}\}$ and $\{u_2^{[n]}\}$ and
  $I'(u_1^{[n]},u_2^{[n]})\to 0$ as $n\to\infty$ imply that
  \begin{eqnarray}\label{dd9}
  |\langle I'(u_1^{[n]},u_2^{[n]}),
  (u_1^{[n]}-u_1,u_2^{[n]}-u_2)\rangle|\le
  \|I'(u_1^{[n]},u_2^{[n]})\|\|(u_1^{[n]}-u_1,u_2^{[n]}-u_2)\|\to
  0,\mbox{ as } n\to\infty.
  \end{eqnarray}
  Moreover, by (G1) and H\"older inequality, we have
 \begin{eqnarray}\label{dd2}
 &     &  \left|\int_\Omega G_{z_1}(x,u_1^{[n]},u_2^{[n]})(u_1^{[n]}-u_1)dx\right|\nonumber\\
 & \le & \int_\Omega c_1(1+|u_1^{[n]}|^{s_1})|u_1^{[n]}-u_1|dx\nonumber\\
 & \le & \left(\int_\Omega c_1^{p'}(1+|u_1^{[n]}|^{s_1})^{p'}ds\right)^{1/p'}\left(\int_\Omega |u_1^{[n]}-u_1|^pdx\right)^{1/p}\nonumber\\
 & \le & 2^{\frac{p'-1}{p'}}c_1\left(\int_\Omega (1+|u_1^{[n]}|^{p's_1})dx\right)^{1/p'}\left(\int_\Omega |u_1^{[n]}-u_1|^pdx\right)^{1/p}.
  \end{eqnarray}
 Since $p's_1\le \frac{pN}{N-p}$, where $p'>1$ satisfies $\frac{1}{p'}+\frac{1}{p}=1$, then Lemma 2.1 and the boundedness of  $\{u_1^{[n]}\}$  imply
 that $\int_\Omega|u_1^{[n]}|^{p's_1}dx$ is bounded. Thus by (\ref{dd4}) and
 (\ref{dd2}), we obtain that
 \begin{eqnarray}\label{dd7}
        \int_\Omega G_{z_1}(x,u_1^{[n]},u_2^{[n]})(u_1^{[n]}-u_1)dx\to  0, \quad \mbox{as } n \to\infty.
  \end{eqnarray}
 Similarly, noting that $q's_2\le \frac{qN}{N-q}$, where $q'>1$ satisfies $\frac{1}{q'}+\frac{1}{q}=1$,   we obtain that
  \begin{eqnarray}\label{dd6}
        \int_\Omega G_{z_2}(x,u_1^{[n]},u_2^{[n]})(u_2^{[n]}-u_2)dx\to  0, \quad \mbox{as }n \to\infty.
  \end{eqnarray}
  Note that $a_1\in C(\bar{\Omega},\mathbb R)$. Then by
  ($\mathcal{M}$) and H\"older inequality, we have
  \begin{eqnarray*}
   &     &  \left|\Big[M_3\left(\int_\Omega a_1(x)|u_1^{[n]}|^p dx\right)\Big]^{p-1}\int_\Omega a_1(x)|u_1^{[n]}|^{p-2}(u_1^{[n]},u_1^{[n]}-u_1)dx \right|\nonumber\\
   & \le &  (m^*)^{p-1}\max_{x\in\bar{\Omega}} a_1(x) \int_\Omega   |u_1^{[n]}|^{p-1}|u_1^{[n]}-u_1|dx\nonumber\\
    & \le &  (m^*)^{p-1}\max_{x\in\bar{\Omega}} a_1(x) \left(\int_\Omega   |u_1^{[n]}|^{p}dx\right)^{1/p'}\left(\int_\Omega    |u_1^{[n]}-u_1|^pdx\right)^{1/p}.
  \end{eqnarray*}
  Then Lemma 2.1 and (\ref{dd4}) imply that
 \begin{eqnarray}\label{dd10}
    \left|\Big[M_3\left(\int_\Omega a_1(x)|u_1^{[n]}|^p dx\right)\Big]^{p-1}\int_\Omega a_1(x)|u_1^{[n]}|^{p-2}(u_1^{[n]},u_1^{[n]}-u_1)dx
    \right|\to 0,\quad\mbox{as }n\to\infty.
  \end{eqnarray}
  Similarly, we also have
  \begin{eqnarray}\label{dd11}
    \left|\Big[M_4\left(\int_\Omega a_2(x)|u_2^{[n]}|^q dx\right)\Big]^{q-1}\int_\Omega a_2(x)|u_2^{[n]}|^{q-2}(u_2^{[n]},u_2^{[n]}-u_2)dx
    \right|\to 0,\quad\mbox{as }n\to\infty.
  \end{eqnarray}
   Thus, (\ref{dd8}), (\ref{dd9}), (\ref{dd7}), (\ref{dd6}), (\ref{dd10}) and (\ref{dd11}) imply
 that
  \begin{eqnarray}\label{dd12}
    &    &  \Big[M_1\left(\int_\Omega|\nabla u_1^{[n]}|^p dx\right)\Big]^{p-1}\int_\Omega|\nabla u_1^{[n]}|^{p-2}(\nabla    u_1^{[n]},\nabla  u_1^{[n]}-\nabla u_1)dx\nonumber\\
    &    &    +\Big[M_2\left(\int_\Omega|\nabla u_2^{[n]}|^q dx\right)\Big]^{q-1}\int_\Omega|\nabla u_2^{[n]}|^{q-2}(\nabla    u_2^{[n]},\nabla u_2^{[n]}-\nabla u_2)dx \nonumber\\
    &    & \to 0, \mbox{ as }n\to \infty.
  \end{eqnarray}
 Then by ($\mathcal{M}$), we have
 \begin{eqnarray}\label{dd13}
    &     &  \frac{1}{\max\{(m^*)^{p-1},(m^*)^{q-1}\}}\Big[M_1\left(\int_\Omega|\nabla u_1^{[n]}|^p dx\right)\Big]^{p-1}\int_\Omega|\nabla u_1^{[n]}|^{p-2}(\nabla    u_1^{[n]},\nabla  u_1^{[n]}-\nabla u_1)dx\nonumber\\
    &     &   +\frac{1}{\max\{(m^*)^{p-1},(m^*)^{q-1}\}} \Big[M_2\left(\int_\Omega|\nabla u_2^{[n]}|^q dx\right)\Big]^{q-1}\int_\Omega|\nabla u_2^{[n]}|^{q-2}(\nabla    u_2^{[n]},\nabla u_2^{[n]}-\nabla u_2)dx \nonumber\\
    & \le &  \int_\Omega|\nabla u_1^{[n]}|^{p-2}(\nabla    u_1^{[n]},\nabla  u_1^{[n]}-\nabla u_1)dx+\int_\Omega|\nabla u_2^{[n]}|^{q-2}(\nabla    u_2^{[n]},\nabla u_2^{[n]}-\nabla u_2)dx \nonumber\\
    & \le &  \frac{1}{\min\{(m_*)^{p-1},(m_*)^{q-1}\}}\Big[M_1\left(\int_\Omega|\nabla u_1^{[n]}|^p dx\right)\Big]^{p-1}\int_\Omega|\nabla u_1^{[n]}|^{p-2}(\nabla    u_1^{[n]},\nabla  u_1^{[n]}-\nabla u_1)dx\nonumber\\
    &     &   +\frac{1}{\min\{(m_*)^{p-1},(m_*)^{q-1}\}} \Big[M_2\left(\int_\Omega|\nabla u_2^{[n]}|^q dx\right)\Big]^{q-1}\int_\Omega|\nabla u_2^{[n]}|^{q-2}(\nabla    u_2^{[n]},\nabla u_2^{[n]}-\nabla u_2)dx. \nonumber\\
  \end{eqnarray}
 Hence, (\ref{dd12}) and (\ref{dd13}) imply that
  \begin{eqnarray}\label{dd14}
     \int_\Omega|\nabla u_1^{[n]}|^{p-2}(\nabla    u_1^{[n]},\nabla  u_1^{[n]}-\nabla u_1)dx+\int_\Omega|\nabla u_2^{[n]}|^{q-2}(\nabla    u_2^{[n]},\nabla u_2^{[n]}-\nabla u_2)dx     \to 0, \mbox{ as } n\to \infty.
  \end{eqnarray}
 Moreover, by H\"older inequality, we have
  \begin{eqnarray}\label{dd15}
  &     &   \left|\int_\Omega| u_1^{[n]}|^{p-2}(u_1^{[n]}, u_1^{[n]}-u_1)dx+\int_\Omega|u_2^{[n]}|^{q-2}(u_2^{[n]},u_2^{[n]}-u_2)dx\right|\nonumber\\
  & \le &   \int_\Omega| u_1^{[n]}|^{p-1}| u_1^{[n]}-u_1|dx+\int_\Omega|u_2^{[n]}|^{q-1}||u_2^{[n]}-u_2|dx\nonumber\\
  & \le &   \left(\int_\Omega| u_1^{[n]}|^{p}dx\right)^{\frac{1}{p'}}\left(\int_\Omega| u_1^{[n]}-u_1|^pdx\right)^{\frac{1}{p}}+  \left(\int_\Omega|u_2^{[n]}|^{q}dx\right)^{\frac{1}{q'}}\left(\int_\Omega|u_2^{[n]}-u_2|^qdx\right)^{\frac{1}{q}}.
  \end{eqnarray}
 Then (\ref{dd4}) implies that
 \begin{eqnarray}\label{dd16}
   \int_\Omega| u_1^{[n]}|^{p-2}(u_1^{[n]}, u_1^{[n]}-u_1)dx+\int_\Omega|u_2^{[n]}|^{q-2}(u_2^{[n]},u_2^{[n]}-u_2)dx \to 0, \mbox{ as }n\to  \infty.
  \end{eqnarray}
 Set
 \begin{eqnarray}\label{dd17}
  J(u) &  =  & \frac{1}{p}\left(\int_\Omega |u_1|^pdx+\int_\Omega |\nabla u_1|^pdx\right)+\frac{1}{q}\left(\int_\Omega |u_2|^qdx+\int_\Omega |\nabla  u_2|^qdx\right)\nonumber\\
       &  =  & J_1(u_1)+J_2(u_2),
 \end{eqnarray}
 where
  \begin{eqnarray*}
  &  &  J_1(u_1)=  \frac{1}{p}\left(\int_\Omega |u_1|^pdx+\int_\Omega |\nabla   u_1|^pdx\right),\quad \forall u_1\in W_0^{1,p},\\
  &  &  J_2(u_2)=  \frac{1}{q}\left(\int_\Omega |u_2|^qdx+\int_\Omega |\nabla  u_2|^qdx\right),\quad \forall u_2\in W_0^{1,q}.
  \end{eqnarray*}
 Then we have
 \begin{eqnarray}\label{dd18}
            \langle J'(u^{[n]}),u^{[n]}-u\rangle
 &   =   &    \langle J'(u_1^{[n]},u_2^{[n]}),(u_1^{[n]}-u_1,u_2^{[n]}-u_2)\rangle\nonumber\\
 &   =   &  \int_\Omega (|u_1^{[n]}|^{p-2}u_1^{[n]},u_1^{[n]}-u_1)dx+\int_\Omega (|\nabla u_1^{[n]}|^{p-2}\nabla u_1^{[n]},\nabla u_1^{[n]}-\nabla u_1)dx\nonumber\\
 &       &  +\int_\Omega (|u_2^{[n]}|^{q-2}u_2^{[n]}, u_2^{[n]}-u_2) dx+\int_\Omega  (|\nabla u_2^{[n]}|^{q-2}\nabla u_2^{[n]},\nabla u_2^{[n]}-\nabla u_2)dx\nonumber\\
 &   =   &  \langle J_1'(u_1^{[n]}),u_1^{[n]}-u_1\rangle+\langle J_2'(u_2^{[n]}),u_2^{[n]}-u_2\rangle.
 \end{eqnarray}
 and
  \begin{eqnarray}\label{dd19}
            \langle J'(u),u^{[n]}-u\rangle
 &   =   &    \langle J'(u_1,u_2),(u_1^{[n]}-u_1,u_2^{[n]}-u_2)\rangle\nonumber\\
 &   =   &  \int_\Omega (|u_1|^{p-2}u_1,u_1^{[n]}-u_1)dx+\int_\Omega (|\nabla u_1|^{p-2}\nabla u_1,\nabla u_1^{[n]}-\nabla u_1)dx\nonumber\\
 &       &  +\int_\Omega (|u_2|^{q-2}u_2, u_2^{[n]}-u_2) dx+\int_\Omega  (|\nabla u_2|^{q-2}\nabla u_2,\nabla u_2^{[n]}-\nabla u_2)dx\nonumber\\
 &   =   &  \langle J_1'(u_1),u_1^{[n]}-u_1\rangle+\langle J_2'(u_2),u_2^{[n]}-u_2\rangle.
 \end{eqnarray}
 (\ref{dd14}) and (\ref{dd16}) imply that
  \begin{eqnarray}\label{dd20}
    \langle J'(u^{[n]}),u^{[n]}-u\rangle  = \langle J_1'(u_1^{[n]}),u_1^{[n]}-u_1\rangle+\langle J_2'(u_2^{[n]}),u_2^{[n]}-u_2\rangle \to 0, \mbox{ as }n\to \infty
  \end{eqnarray}
 and (\ref{dd21}) implies that
 \begin{eqnarray}\label{dd23}
    \langle J'(u),u^{[n]}-u\rangle = \langle J_1'(u_1),u_1^{[n]}-u_1\rangle+\langle J_2'(u_2),u_2^{[n]}-u_2\rangle\to 0, \mbox{ as }n\to \infty.
  \end{eqnarray}
  By Lemma 2.3 in \cite{AL}, we know that
 \begin{eqnarray}\label{dd24}
   \langle J_1'(u_1^{[n]})-J_1'(u_1),u_1^{[n]}-u_1\rangle \ge (\|u_1^{[n]}\|_{1,p}^{p-1}-\|u_1\|_{1,p}^{p-1}) (\|u_1^{[n]}\|_{1,p}-\|u_1\|_{1,p})
  \end{eqnarray}
  and
 \begin{eqnarray}\label{dd25}
 \langle J_2'(u_2^{[n]})-J_2'(u_2),u_2^{[n]}-u_2\rangle \ge (\|u_2^{[n]}\|_{1,q}^{q-1}-\|u_2\|_{1,q}^{q-1}) (\|u_2^{[n]}\|_{1,q}-\|u_2\|_{1,q}).
  \end{eqnarray}
 Then it follows from (\ref{dd24}) and (\ref{dd25})  that
 \begin{eqnarray*}
              0
   &  \le  & (\|u_1^{[n]}\|_{1,p}^{p-1}-\|u_1\|_{1,p}^{p-1}) (\|u_1^{[n]}\|_{1,p}-\|u_1\|_{1,p})+ (\|u_2^{[n]}\|_{1,q}^{q-1}-\|u_2\|_{1,q}^{q-1}) (\|u_2^{[n]}\|_{1,q}-\|u_2\|_{1,q})\nonumber\\
   &  \le  &   \langle J'(u^{[n]})-J'(u),u^{[n]}-u\rangle,
  \end{eqnarray*}
  which, together with (\ref{dd20}) and (\ref{dd23}), yields $\|u_1^{[n]}\|_{1,p} \to
  \|u_1\|_{1,p}$ and  $\|u_2^{[n]}\|_{1,q} \to  \|u_2\|_{1,q}$ as $n \to
  \infty$. Note that $W_0^{1,p}$ and $W_0^{1,q}$ are reflexive and  uniformly
  convex (see \cite{Adams}). Then it follows from Kadec-Klee property (see \cite{HS}) that $u_1^{[n]}\to
  u_1$ strongly in  $W_0^{1,p}$  and $u_2^{[n]}\to
  u_2$ strongly in  $W_0^{1,q}$  and so $u^{[n]}\to
  u$ strongly in  $\mathcal{W}$. Thus we have verified that $I$ satisfies the
  (C)-condition. This completes the proof.

 \vskip2mm
 \noindent
 {\bf Lemma 3.5.} {\it Assume that (G2)$'$ and (G4) hold. Then $I$ satisfies the  (C)-condition. }
 \vskip2mm
 \noindent
 {\bf Proof.}   The proof is the same as that of Lemma 3.4 by replacing
 (\ref{b4}) with (\ref{b41}).

 \vskip2mm
 \noindent
 {\bf Proof of Theorem 1.1.} Obviously, by (G0), we have $I(0)=0$ and  $I$ is
 even. Let $E=\mathcal{W}$, $E_1=E_{p,1}^{(1)}\times E_{q,1}^{(1)}$ and $E_2=E_{p,1}^{(2)}\times E_{q,1}^{(2)}$. Then $E_1$ is finite
 dimensional. Combining Lemma 3.1, Lemma 3.3, Lemma 3.4  with Lemma 2.3, we obtain that system (\ref{a1}) has
 infinitely many nontrivial solutions.

 \vskip2mm
 \noindent
 {\bf Proof of Theorem 1.2.} The proof is the same as that of Theorem 1.1 by replacing  Lemma 3.2
 with Lemma 3.1 and replacing  Lemma 3.4 with Lemma 3.5.

  {}

\end{document}